\documentclass{elsart}

\usepackage{amssymb}
\newtheorem{theorem}{Theorem}[section]
\newtheorem{lemma}{Lemma}[section]
\newtheorem{remark}{Remark}[section]

\begin{document}

\begin{frontmatter}



\thanks{Supported by  DFG 436 RUS
113/823/0-1.}

\title{Blow up of smooth highly decreasing at infinity solutions to the
compressible Navier-Stokes equations}


\author{ Olga Rozanova}

\address{Department of Differential Equations \& Mechanics and Mathematics Faculty,
Moscow State University, Moscow, 119992,
 Russia}
\ead{rozanova@mech.math.msu.su}
\begin{abstract}We prove that the smooth solutions to the Cauchy problem
for the  Navier-Stokes equations with  conserved total mass, finite
total energy and finite momentum of inertia  lose the initial
smoothness within a finite time in the case of space of dimension 3
or greater even if the initial data are not compactly supported. The
cases of isentropic and incompressible fluids are also considered.

\end{abstract}

\end{frontmatter}

\section{System, known results and main problem}
\label{1}











The motion of compressible viscous, heat-conductive, 
Newtonian polytropic fluid in ${\mathbb R}\times{\mathbb R}^n,\,
n\ge 1,$ is governed by the compressible Navier-Stokes (NS)
equations
$$\partial_t \rho+{\rm div}_x (\rho u)=0,\eqno(1.1)$$
$$\partial_t(\rho u)+{\rm div}_x (\rho u \otimes u)+\nabla_x p={\rm Div} T,\eqno(1.2)$$
$$\partial_t\left(\frac{1}{2}\rho |u|^2+\rho e\right)+{\rm div}_x \left((\frac{1}{2}\rho |u|^2+\rho e+p)u
\right)={\rm div}(Tu)+k\Delta_x \theta, \eqno(1.3)     $$ where
$\rho, \, u=(u_1,...,u_n),\, p, \, e,\,\theta \,$ denote the
density, velocity, pressure, internal energy and absolute
temperature, respectively, $\,T\,$ is the stress tensor given by the
Newton law
$$T=T_{ij}=\mu \,(\partial_iu_j+\partial_j u_i)+\lambda \,{\rm div} u\, \delta_{ij},\eqno (1.4)$$
where the constants $\mu$ and $\lambda$ are the coefficient of
viscosity and the second coefficient of viscosity, $k\ge 0$ is the
coefficient of heat conduction. We denote ${\rm Div}$ and $\rm div$
the divergency of tensor and vector, respectively. We assume that
$\mu > 0, \,\lambda+\frac{2}{n}\mu >0.$

The state equations have the forms
$$p=R\rho \theta,\quad e=c\theta,\quad
p=A\exp\left(\frac{S}{c}\right)\rho^\gamma.\eqno(1.5)
$$
Here $A>0$ is a constant, $R$ is the universal gas constant,
 $S=\log e - (\gamma-1)\log
\rho$ is the specific entropy, $c=\frac{R}{\gamma-1},$ $\gamma>1$ is
the specific heat ratio,

The state equations (1.5) imply
$$\,p=(\gamma-1)\rho e,\eqno(1.6)$$
which allows us to consider (NS) as a system for the unknown
$\rho,\,u,\,p.$ Indeed,
 from (NS) and (1.6) it follows that
$$\partial_t p +
(u,\nabla_x p)+\gamma p\,{\rm div}u=(\gamma-1)\,\sum\limits_{i,j=1}^n\,T_{ij}
\partial_j u_i+\frac{k}{R}\,\Delta\frac{p}{\rho}.\eqno(1.7)$$
Thus, therefore further we shall consider the system (1.1, 1.2,
1.7), denoted (NS*) for short.



(NS*) is supplemented with the initial data
$$(\rho,u,p)\Bigl|_{t=0}=(\rho_0(x),u_0(x),p_0(x))\in H^m({\mathbb
R}^n),\,m>[n/2]+2.\eqno(1.8)$$


We also consider  the isentropic case where the fluid obeys
equations (1.1), (1.2) and $p=A \rho ^\gamma,$ we call this system
(NSI) for short.

In the absence of vacuum, the local existence of classical solutions
is known. Namely, in \cite{Nash} it is proved that there exist
classical solutions, having the H\"older continuous second
derivatives with respect to space variables and the first ones with
respect to time. In \cite{VolpertKhudiaev} the system of equations
of viscous compressible fluid is considered as a particular case of
combined systems of differential equations. The consideration is
performed in the Sobolev spaces $H^m$ with a sufficiently large $m.$
The uniqueness of the solution was proved earlier in \cite{Serrin}.
The existence and uniqueness of local strong solutions in the case
where the initial density need not be positive and may vanish in an
open set were proved recently in\cite{ChoKim}.

At the same time there exists a major open problem: to prove or
disprove that a smooth solution to (NS) (in higher space dimensions)
exists globally in time.  There are partial results concerning the
Cauchy problem for (NS) away from a vacuum. In
\cite{MatsumuraNishida} it is proved that if there exists a constant
 $\bar \rho=const >0 $ such that $(\rho_0-\bar \rho, u_0,S_0)\in H^m({\mathbb R}^n)\cap L^1({\mathbb
R}^n),\,$ and the norm $\,\|\rho_0-\bar\rho,
u_0,S_0\|_{H^m},\,m>[n/2]+2,\,$ is suitably small, then
 global solution to (NS) from $C^1([0,\infty),
H^m({\mathbb R}^n))$ exists. In \cite{AntontsevKazhikhovMonakhov} it
is shown that for $n=1$ the global existence takes place without
assumptions on the smallness of the norm $\|\rho_0-\bar\rho,
u_0,S_0\|_{H^m}.$

However if the initial density $\rho_0$ is compact, then in
arbitrary space dimensions no solution to (NS) from $C^1([0,\infty),
H^m({\mathbb R}^n))$ exists (\cite{Xin}). This blowup result depends
crucially on the assumption about  compactness of support of the
initial density. It does not seem to solve in a negative way the
question of regularity for (NS). Indeed, (NS) is a model of
non-dilute fluids where the density is bounded below away from zero,
and therefore it is natural to expect the problem to be ill-posed
when vacuum regions are present at the initial time. At the same
time, the conservation of mass in the whole space requires a
decrease of the density down to zero.

In \cite{Xin} the author notes that the global smooth solution to
(NS) seems to exist at least for small data in the case where
initial vacuum appears only at infinity.

Nevertheless, in \cite{ChoJa} a sufficient condition for the blow-up
in case that the initial density is positive but has a decay at
infinity was found. In this work it was expected a specific time
decay of the velocity component that seems to be reasonable for a
density away from zero.  Further, in \cite{Ro} it was proved that if
the solution with finite moment of inertia to the
(non-heat-conductive) Navier-Stokes system is smooth globally in
time, then the solution components grow as $t\to\infty$ at least as
a certain function specific for every component $\phi(t,x),$ where
$\phi$ stands either for the density, pressure, velocity or gradient
of velocity. Namely, we observe all trajectories of particles $x(t)$
that leave in a finite time a ball of finite radius $R_0$ and find
that provided the solution is smooth, the inequality $\,|\phi
(t,x(t))|\le M(t) |x(t)|^\lambda, \,\,t\to\infty,\,$ holds, with a
continuous function $M(t)$ and a constant $\lambda, \,$ specific for
every $\phi.$  Both in \cite{ChoJa} and \cite{Ro} the answer whether
the solution blows up or not does not depend on the initial data,
but only on the prescribed decay at infinity.

\vskip0.5cm

Thus, the question remains: is it true that in the case where the
support of all components of initial data coincides with the whole
space the global smooth solution exists for any smooth initial data?

Below we are going to show that, generally speaking, the answer is
negative.


\vskip1cm

\section{Integral functionals, solution with decreasing components  and the statement of main theorem}

\vskip1cm

The system (NS) is the differential form of conservation laws for
the material volume $\Omega(t)$; it expresses  conservation of mass
$$m=\int\limits_{\Omega(t)}\rho \, d x,$$  and balance of momentum
$$P=\int\limits_{\Omega(t)}\rho u\, d x,$$ and  total energy
$$\mathcal E=\int\limits_{\Omega(t)}\left(\frac{1}{2}\rho |u|^2+\rho
e\right)\, d x \,=E_{k}(t)+E_{i}(t).$$ Here $E_k(t)$ and $E_i(t)$
are the kinetic and internal components of energy, respectively.

 If we regard $\Omega (t)={\mathbb R}^n,$ the conservation of
mass,  momentum and energy takes place provided the components of
the solution  decrease at infinity sufficiently quickly.

\vskip1cm

 Let us introduce the  functionals
$$G(t)=\frac{1}{2} \int\limits_{{\mathbb
R}^n}\rho(t,x)|{ x}|^2\,dx,\qquad F(t)=\int\limits_{{\mathbb R}^n}
({u},{ x})\rho\,dx,$$ where the first one is the momentum of
inertia, the scalar product of vectors is denoted as $(.,.)$ .

For technical reasons we impose the following conditions of decay on
the solution  components to (NS*) as $|{ x}|\to \infty$ at every
fixed $t\in {\mathbb R}_+$:
$$
\rho = O\left(\frac{1}{|x|^{n+2+\varepsilon}}\right),\qquad p =
O\left(\frac{1}{|x|^{n+\varepsilon}}\right),\quad
\varepsilon>0,\eqno(2.1)$$
$$|{u
}|=o\left(\frac{1}{|x|^{n-1}}\right), \qquad |D{
u}|=o\left(\frac{1}{|x|^{n}}\right).\eqno(2.2)
$$
If $k\ne 0,$ we require additionally
$$
|D\theta|=o\left(\frac{1}{|x|^{n-1}}\right),\quad |{ x}|\to
\infty,\quad t\in {\mathbb R}_+.\eqno(2.3)
$$

 One can easily verify that
these requirements guarantee the conservation of the mass $m,$
energy $\mathcal E,$ momentum $P$  on solutions to (NS*) and ensure
a convergence of the momentum of mass $G(t).$

We impose no restriction on the solution support, however the decay
of the solution as $|{ x}|\to \infty $ in the class  considered is
greater than it is necessary  for belonging to $C^1([0,T),
H^m({\mathbb R}^n)).$

\vskip1cm

\begin{defn}
We will say that a solution $(\rho, u,p)$ to the Cauchy problem
(1.1, 1.2, 1.7), (1.8) belongs to the class $\mathfrak K$ if it has
the following properties for all $t\ge 0$ :

(i)\,\,\quad the solution is classical;

(ii)\,\quad the solution decays at infinity according to (2.1 --
2.3);

(iii)\quad  $\rho(t,x)\ge 0;\, $

(iv) \quad  $\displaystyle \frac{d S(t,x)}{dt}=\sigma(t,x)\ge 0,\,$
\quad $\big\|\sigma(t,x)\big\|_{L^\infty({\mathbb R}^n)}=
o(t^\alpha), \quad t\to\infty,\quad \mbox{where}\quad\alpha =
\frac{(\gamma-1)n^2+n-2}{n}, $ if $\gamma\le 1+\frac{2}{n},$ and
$\alpha=\frac{3n-2}{n},$ otherwise.
\end{defn}

\medskip
\begin{remark}In
particular, $(iv)$ results that $S(t,x)\ge S_0 =\displaystyle
\inf\limits_{x\in {\mathbb R}^n} S(0,x)$. For (NSI)  condition
$(iv)$ holds trivially.
\end{remark}

\vskip1cm

\begin{remark}
 From (1.7) we have
$$\frac{p}{R}\left(\partial_t S +(u,\nabla_x
S)\right)=\sum\limits_{i,j=1}^n\, T_{ij}\partial_j
u_i+k\Delta\theta.$$ Evidently,  if the right-hand side is
non-negative, then $S(t,x)\ge S_0 $ for all $t>0.$
 The first item on the right-hand side is non-negative. Indeed,
according to \cite{Xin}
$$
\sum\limits_{i,j=1}^n\, T_{ij}\partial_j u_i=
\sum\limits_{i,j=1}^n\,\partial_j (T_{ij} u_i)-
\sum\limits_{i,j=1}^n\,  u_i \,\partial_j T_{ij}=$$
$$
2\mu\,\sum\limits_{i,j=1}^n\, (\partial_j u_j)^2+ \lambda ({\rm
div}\,u)^2+\mu\,\sum\limits_{i\ne j}^n\, (\partial_j u_i)^2+
2\mu\,\sum\limits_{i>j}^n\, (\partial_j u_i)(\partial_i
u_j).\eqno(2.4)
$$
For $\lambda\ge 0$ this implies the nonnegativity of the right-hand
side of (2.4), for $\lambda<0$ (2.4) can be estimated from below by
$$
(2\mu+n\lambda)\,\sum\limits_{i,j=1}^n\, (\partial_j u_j)^2+
\mu\,\sum\limits_{i\ne j}^n\, (\partial_j u_i)^2+
2\mu\,\sum\limits_{i>j}^n\, (\partial_j u_i)(\partial_i u_j);
$$
this expression is nonnegative by assumption on the viscosity
coefficients.

 Therefore we can guarantee the
uniform boundedness of entropy from below in the case $k=0$ and in
the isothermal case $\theta=const.$
\end{remark}

\vskip1cm

\begin{theorem} Let $n\ge 3,\,\gamma\ge \frac{2n}{n+2},$ the momentum $P\ne 0.$ If $\displaystyle \inf\limits_{x\in{\mathbb
R}^n} \,S(0,x)>-\infty, $ then there exists no global in time
solution to (NS*) from the class $\mathfrak K.$
\end{theorem}


Analyzing condition $(iv)$ we can see that the nondecreasing of
entropy along trajectories seems natural, whereas the upper bound on
the growth of entropy can be considered as unreasonable. So, we can
re-formulate Theorem 2.1 as follows:

\begin{theorem} Let us assume $n\ge 3,\,\gamma\ge \frac{2n}{n+2}$
and $P\ne 0.$ Then any solution to (NS*) with properties $(i, ii,
iii)$ such that the entropy does not decrease along the particles
trajectories, blows up in a finite or infinite time. If the solution
keeps smoothness for all $t>0,$ then $\|S\|_{L^\infty({\mathbb
R}^n)},$ $\|p\|_{L^\infty({\mathbb R}^n)},$ $\|{\rm div}
u\|_{L^\infty({\mathbb R}^n)}$ rise at least as $O(t^{\alpha+1})$ as
$t\to \infty$  (the constant $\alpha>0$ is indicated in condition
$(iv)$.)
\end{theorem}

For the isentropic case we have the following version of Theorem
2.1.

\begin{theorem} Let $n\ge 3,\,\gamma\ge \frac{2n}{n+2}.$
If initial data $(\rho_0(x), u_0(x))$ satisfy that the momentum
$P\ne 0,$ then the solution to (NSI) from the class $\mathfrak K$
cannot exist for all $t>0.$
\end{theorem}

\vskip1cm

\section{Proof of the theorems}

\vskip1cm We begin with the following  extension of formula obtained
in \cite{Chemin} to the viscous case.
\begin{lemma}
For  solutions to (1.1, 1.2, 1.7) with properties $(i), (ii) $
$$G'(t)=F(t),\eqno(3.1)$$
$$F'(t)=2 E_k(t)+n(\gamma-1) E_i(t).\eqno(3.2)$$
\end{lemma}

\begin{pf}
The lemma can be proved by direct calculation using the general
Stokes formula and taking into account the decay assumption $(ii).$
For example, from (1.1) we get
$$G'(t)=\frac{1}{2}\,\int\limits_{{\mathbb R}^n}\rho'_t|{
x}|^2\,dx=$$$$=\int\limits_{{\mathbb R}^n} ({ u},{
x})\rho\,dx-\lim_{R\to \infty}\int\limits_{S_R}
  {(u,x)}\rho \frac{|{ x}|}{2}\,dS_R=$$
$$=\int\limits_{{\mathbb R}^n}({u},{x})\rho\,dx,$$
where $S_R$ is a sphere of radius $R$ with the center at the origin.
This proves equality (3.1).
\end{pf}

Then we get two-sided estimates of $G(t).$

\begin{lemma}
If $\gamma\le 1+\frac{2}{n}, $ then for solutions with properties
$(i), (ii) $ the estimates
$$\frac{n(\gamma-1)}{2}\mathcal E t^2+F(0)t+G(0)\le G(t)\le  \mathcal E t^2+F(0)t+G(0)\eqno(3.3)$$
(for (NS*)) and
$$\frac{P^2}{2m}t^2+F(0)t+G(0)\le G(t)\le  \mathcal E(0) t^2+F(0)t+G(0)\eqno(3.4)$$
(for (NSI)) hold.

If $\gamma> 1+\frac{2}{n}, $ then we have
$$\mathcal E t^2+F(0)t+G(0)\le G(t)\le \frac{n(\gamma-1)}{2} \mathcal E t^2+F(0)t+G(0)\eqno(3.5)$$
for (NS*) and
$$\frac{P^2}{2m}t^2+F(0)t+G(0)\le G(t)\le \frac{n(\gamma-1)}{2} \mathcal E(0) t^2+F(0)t+G(0)\eqno(3.6)$$
for (NSI).

\end{lemma}
\begin{pf}
First of all (3.2) result
$$
G''(t)=2E_k(t)+n(\gamma-1)E_i(t)=$$
$$=2\mathcal E (t) -
(2-n(\gamma-1))E_i(t)=n(\gamma-1)\mathcal
E(t)+(2-n(\gamma-1))E_k(t).\eqno(3.7)
$$
 In the case of (NS*) the total energy
$\mathcal E$ is constant, therefore (3.3) and (3.5) follow from
(3.7) after integration if we take into account  nonnegativity of
$E_k(t)$ and $E_i(t).$ Estimate (3.3) one can find in \cite{Chemin}
for the zero viscosity.

In the case of (NSI) the total energy is only non-increasing, since
(1.1) and (1.2) for a constant entropy result in
$$\mathcal E'(t)\le -\nu \int \limits_{{\mathbb R}^n}|D u|^2\, dx\le 0,\eqno(3.8)$$
with some positive constant $\nu.$ Then to get the two-sided
estimates of $G(t)$ from (3.7) we can use
 the non-increasing of total energy
$\mathcal E(t)\le \mathcal E(0),$ the nonnegativity of $E_i(t)$ and
 the estimate $E_k(t)\ge\frac{P^2}{2m},$ that follows from the H\"older inequality. This gives (3.4) and (3.6).

\end{pf}

The next step is two-sides estimate of $E_i(t).$

\vskip1cm

\begin{lemma}
 For
solutions from the class $\mathfrak K$ to (NS*)  for sufficiently
large $t$ we have
$$\frac{C_1}{G^{(\gamma-1)n/2}(t)}\le E_i(t) \le \frac{C_2}{G^{(\gamma-1)n/2}(t)},\eqno(3.9)$$
for $\gamma\le 1+\frac{2}{n},$ and
$$\frac{C_1}{G^{(\gamma-1)n/2}}\le E_i(t) \le \frac{C_2}{G(t)},\eqno(3.10)$$
for $\gamma > 1+\frac{2}{n},$ with constants $C_1, C_2 \,(C_1\le
C_2).$ The same estimate holds also for (NSI) provided $P\ne 0.$
\end{lemma}

\begin{pf}
The lower estimate is due to \cite{Chemin}. It follows from the
inequality
$$
\|f\|_{L^1({\mathbb R}^n;\,dx)}\,\le\,C_{\gamma,
n}\|f\|^{\frac{2\gamma}{(n+2)\gamma-n}}_{L^\gamma({\mathbb
R}^n;\,dx)}\,\|f\|^{\frac{n(\gamma-1)}{(n+2)\gamma-n}}_{L^1({\mathbb
R}^n;\,|x|^2\,dx)},
$$
together with the lower estimate of internal energy
$$E_i(t)\ge \frac{A e^{S_0/c}}{\gamma-1}\,\int\limits_{{\mathbb R}^n}\,\rho^\gamma(t,x)\,dx.$$
The latter inequality follows from the state equations in (1.5). The
constant
$$C_1=A \frac{e^{S_0/c}}{\gamma-1}(m C_{\gamma,
n}^{-1})^\frac{\gamma(n+2)-n}{2}$$ with
 $$C_{\gamma,
n}=\left(\frac{2\gamma}{n(\gamma-1)}\right)^
{\frac{n(\gamma-1)}{(n+2)\gamma-n}} +
\left(\frac{2\gamma}{n(\gamma-1)}\right)^
{\frac{-2\gamma}{(n+2)\gamma-n}}.$$

The method of the upper estimate of $E_i(t)$ is also similar to
\cite{Chemin}. Namely, let us consider the function
$Q(t)=4G(t)\mathcal E(t)-F^2(t).$ The H\"older inequality gives
$F^2\le 4 G(t) E_k(t),$ therefore $\mathcal E(t)=E_k(t)+E_i(t)\ge
E_i(t)+\frac{F^2(t)}{4G(t)}$ and
$$
E_i(t)\le \frac {Q(t)}{4 G(t)}.\eqno(3.11)
$$
We notice also that $Q(t)>0$ provided the pressure does not equal to
zero identically. Then taking into account (3.1), (3.2) and (3.7) we
have
$$
Q'(t)=4 G'(t) \mathcal E(t)-2 G'(t) G''(t)+ 4G(t)\mathcal
E'(t)=$$$$= 2(2-n(\gamma-1))\,G'(t)E_i(t)+4G(t)\mathcal
E'(t).\eqno(3.12)
$$
Further, one can see from (3.2) that in the (NS*) case $G'(t)>0$
beginning from a positive $t_0$ for all initial data, whereas for
(NSI) equality (3.2) result in
$$G''(t)\ge \frac{P^2}{m},$$
and we can guarantee the positivity of $G'(t)$ for sufficiently
large $t$ only for $P\ne 0.$ Thus, for $\gamma\le 1+\frac{2}{n}$
(3.11), (3.12) (and (3.8) for (NSI)) result
$$\frac{Q'(t)}{Q(t)}\le \frac{2-n(\gamma-1)}{2}\,\frac{G'(t)}{G(t)}.\eqno(3.13)$$
Then (3.13) and (3.11) give
$$E_i(t)\le \frac{C_2}{G^{(\gamma-1)n/2}(t)}, \quad C_2=\frac{Q(0)G^{(\gamma-1)n/2}(0)}{4}.\eqno(3.14)$$
If $\gamma>1+\frac{2}{n},$ then from (3.12) taking into account the
lower estimate of $E_i(t)$ we get
$$
Q'(t)\le -2 ((\gamma-1)n-2)\,\frac{C_1}{G^{(\gamma-1)n/2}(t)},
$$
$$
Q(t)\le \tilde C+{4 \, C_1}\,{G^{1-(\gamma-1)n/2}(t)},\quad \tilde
C=Q(0)-{4 \, C_1}\,{G^{1-(\gamma-1)n/2}(0)},
$$
and, at last, from (3.11)
$$E_i(t)\le \frac{\tilde C}{4\,G(t)}+{ C_1}\,{G^{-(\gamma-1)n/2}(t)}.$$
Thus, since the leading term in the right hand side of the latter
inequality is the first one, beginning from a moment $t>0$ we get
$$E_i(t)\le \frac{C_2}{G(t)}, \quad  C_2= \frac{\tilde C}{4}+C_1.$$
The proof is over.

\end{pf}

The next lemma is a key point of the theorem's proof.

\begin{lemma} Let $n\ge 3,\,\gamma\ge \frac{2n}{n+2}.$ If $|P|\ne 0,$ then there
exists a positive constant $K$ such that for the solutions of the
class $\mathfrak K$ the following inequality holds:
$$
\int\limits_{{\mathbb R}^n}\,|Du|^2\,dx \ge
\,K\,E_i^{-\frac{n-2}{n(\gamma-1)}}(t). \eqno(3.15)
$$

\end{lemma}

\begin{pf}
First of all we use the H\"older inequality to get
$$
|P|=\Big|\int\limits_{{\mathbb R}^n}\,\rho u \,dx\Big|\le
\left(\int\limits_{{\mathbb R}^n}\,\rho^{\frac{2n}{n+2}}
\,dx\right)^{\frac{n+2}{2n}}\left(\int\limits_{{\mathbb R}^n}\,
|u|^{\frac{2n}{n-2}} \,dx\right)^{\frac{n-2}{2n}}\le
$$
$$\le
e^{-\frac{n-2}{2n(\gamma-1)}\,\frac{S_0}{c}}\,\left(\int\limits_{{\mathbb
R}^n}\,e^{\frac{n-2}{(n+2)(\gamma-1)}\,\frac{S}{c}}\rho^{\frac{2n}{n+2}}
\,dx\right)^{\frac{n+2}{2n}}\left(\int\limits_{{\mathbb R}^n}\,
|u|^{\frac{2n}{n-2}} \,dx\right)^{\frac{n-2}{2n}}.\eqno(3.16)
$$
Further, using the Jensen inequality we have for
$\frac{(\gamma-1)(n+2)}{n-2}\ge 1$ (or $\gamma\ge\frac{2n}{n+2}$)
$$
\left(\frac{1}{m}\,\int\limits_{{\mathbb
R}^n}\,e^{\frac{n-2}{(n+2)(\gamma-1)}\,\frac{S}{c}}\,\rho^{\frac{2n}{n+2}}
\,dx\right)^{\frac{(\gamma-1)(n+2)}{n-2}}\le
\frac{\int\limits_{{\mathbb R}^n}e^{\frac{S}{c}}\rho^\gamma
\,dx}{m}=\frac{(\gamma-1) E_i(t)}{mA}.
$$
Thus,
 the latter inequality and (3.16) give
$$
|P|\le K_1\, \left(E_i(t)\right)^{\frac{n-2}{2n(\gamma-1)}}
\left(\int\limits_{{\mathbb R}^n}\, |u|^{\frac{2n}{n-2}}
\,dx\right)^{\frac{n-2}{2n}},\eqno(3.17)
$$
with the positive constant $K_1$ that depends on $\gamma, n, m,
S_0.$ Further, we take into account of the inequality
$$
\left(\int\limits_{{\mathbb R}^n}\,
|u|^{\frac{2n}{n-2}}\,dx\right)^{\frac{n-2}{n}} \le K_2
\,\int\limits_{{\mathbb R}^n}\, |Du|^2 \,dx,\eqno(3.18)
$$
where the constant $K_2>0 $ depends on $n,$  $\,n\ge 3.$ The latter
inequality holds for $u\in H^1({\mathbb R}^n)$ (\cite{Hebey}, p.22)
and follows from the Sobolev embedding.

Thus, the lemma statement, the inequality (3.15), follows from
(3.17), (3.18), with the constant $K=\frac{|P|^2 }{K_1^2\,K_2}.$

\end{pf}

\vskip1cm

We begin from the isentropic case (NSI).

{\it Proof of theorem 2.3.}
From (3.8), (3.15) we have
$$\mathcal E'(t)\le -\nu \,K \,(E_i(t))^{-\frac{n-2}{n(\gamma-1)}}.\eqno (3.19)  $$
Together with (3.4), (3.6) and Lemma 3.3 inequality (3.19) implies
$$
\mathcal E'(t)\le-\nu\, K\,
C_2^{\frac{2-n}{n(\gamma-1)}}\,G^{\frac{n-2}{2n}}(t)\le - L
\,t^{\frac{n-2}{n}},
$$

with a positive constant $L,$ for $\gamma\le 1+\frac{2}{n},$ and
$$
\mathcal E'(t)\le-\nu\, K\,
C_2^{\frac{2-n}{n(\gamma-1)}}\,G^{\frac{n-2}{n(\gamma-1)}}(t)\le - L
\, t^{\frac{2(n-2)}{n(\gamma-1)}},
$$
for $\gamma> 1+\frac{2}{n}.$ In both cases this contradicts to the
non-negativity of $\mathcal E(t).$ Thus, theorem 2.3 is proved.

\medskip

{\it Proof of theorem 2.1.} Let us remind that in the case of (NS*)
the total energy $\mathcal E$ is constant for the solutions of the
class $\mathfrak K.$  However, the derivatives of both kinetic and
internal components of the total energy can be estimated. Namely,
taking into account the state equation in (1.5) we have
$$
\frac{d \,E_k(t)}{dt}=\int\limits_{\mathbb{R}^n} \,(u, {\rm Div}
T)\,dx -\int\limits_{\mathbb{R}^n} \,(u, \nabla p)\,dx\le
$$
$$
\le -\nu\,\int\limits_{\mathbb{R}^n} \,|Du|^2\,dx- \frac{A}{c}\,
\int\limits_{\mathbb{R}^n} \,e^{\frac{S}{c}}\,\rho^\gamma\,(u,\nabla
S)\,dx - A\,\int\limits_{\mathbb{R}^n} \,e^{\frac{S}{c}}\,(u,\nabla
\rho^\gamma)\,dx,\eqno(3.20)
$$
with a positive constant $\nu.$ Further, together with (1.1)  we
obtain
$$
\frac{d \,E_i(t)}{dt}=\frac{d}{dt}\left(\int\limits_{\mathbb{R}^n}
\,\frac{A\,e^{\frac{S}{c}}\,\rho^\gamma}{\gamma-1}\,dx \right)=$$$$=
\frac{A}{c(\gamma-1)}\,\int\limits_{\mathbb{R}^n}
\,e^{\frac{S}{c}}\,\rho^\gamma \partial_t S\,dx+
\frac{A\gamma}{\gamma-1}\,\int\limits_{\mathbb{R}^n}
\,e^{\frac{S}{c}}\,\rho^{\gamma-1} \partial_t \rho\,dx=$$$$=
\frac{A}{c(\gamma-1)}\,\int\limits_{\mathbb{R}^n}
\,e^{\frac{S}{c}}\,\rho^\gamma \frac{dS}{dt}\,dx-
\frac{A\gamma}{\gamma-1}\,\int\limits_{\mathbb{R}^n}
\,e^{\frac{S}{c}}\,\rho^{\gamma-1} {\rm div}(\rho u)\,dx
-$$$$-\frac{A}{c(\gamma-1)}\,\int\limits_{\mathbb{R}^n}
\,e^{\frac{S}{c}}\,\rho^\gamma (u,\nabla S)\,dx\le
 \eqno(3.21)$$$$\le \frac{\big\|\sigma(t,x)\big\|_{L^\infty({\mathbb
R}^n)}}{c}\,E_i(t)+
\frac{A\gamma}{\gamma-1}\,\int\limits_{\mathbb{R}^n} \,\rho
u\,\nabla(e^{\frac{S}{c}}\,\rho^{\gamma-1})\,dx-$$$$
-\frac{A}{c(\gamma-1)}\,\int\limits_{\mathbb{R}^n}
\,e^{\frac{S}{c}}\,\rho^\gamma (u,\nabla S)\,dx=
$$
$$=\frac{\big\|\sigma(t,x)\big\|_{L^\infty({\mathbb R}^n)}}{c}\,E_i(t)
+\frac{A}{c}\, \int\limits_{\mathbb{R}^n}
\,e^{\frac{S}{c}}\,\rho^\gamma\,(u,\nabla S)\,dx +
A\,\int\limits_{\mathbb{R}^n} \,e^{\frac{S}{c}}\,(u,\nabla
\rho^\gamma)\,dx.
$$
At last, from (3.20), (3.21), the condition $(iv),$ and estimates
(3.9), (3.10), (3.3), (3.5) one can get for $\gamma\le
1+\frac{2}{n}$
$$
0= \frac{d}{dt}(E_k(t)+E_i(t))\le$$$$ \le
-\nu\,\int\limits_{\mathbb{R}^n} \,|Du|^2\,dx +\frac
{\big\|\sigma(t,x)\big\|_{L^\infty({\mathbb R}^n)}}{c}\,E_i(t)\le$$
$$\le -\nu \,K\,
C_2^{\frac{2-n}{n(\gamma-1)}}\,G^{\frac{n-2}{2n}}(t) + o(t^\alpha)
\,G^{-\frac{n(\gamma-1)}{2}}(t)\le$$
$$\le -L\, t^{\frac{n-2}{n}}+ o(t^{\alpha -
n(\gamma-1)}),
$$
with some positive constant $L.$ Analogously, for $\gamma>
1+\frac{2}{n},$
$$
0=\frac{d}{dt} (E_k(t)+E_i(t))\le  -\nu\, K
\,C_2^{\frac{2-n}{n(\gamma-1)}}\,G^{\frac{n-2}{2n}}(t) +
o(t^\alpha)\, G^{-1}(t)\le -L\, t^{\frac{n-2}{n}}+ o(t^{\alpha -
2}).
$$

If we take $\alpha$ from the condition $(iv),$ we get a
contradiction that proves the theorem.

\medskip

{\it Proof of theorem 2.2.} It remains to prove the upper estimates
for $\|p\|_{L^\infty({\mathbb R}^n)}$ and $\|{\rm div}
u\|_{L^\infty({\mathbb R}^n)}.$ It is easy to compute that
$$E_i'(t)=-\int\limits_{{\mathbb R}^n} \,p\,\,{\rm div} u\,dx + (\mu +\lambda)\,
\int\limits_{{\mathbb R}^n} |{\rm div} u|^2\,dx + \mu
\int\limits_{{\mathbb R}^n} \,|D u|^2\,dx.\eqno (3.22)$$ Let us
denote $f_1(t)=\|{\rm div} u\|_{L^\infty({\mathbb R}^n)}$ and
$f_2(t)=\|p\|_{L^\infty({\mathbb R}^n)}.$ As in the proof of Theorem
2.1 taking into account Lemma 3.3 we get from (3.22)
$$
E_i'(t)\ge -(\gamma-1)\, E_i(t) \, f_1(t)+ \mu
\,\int\limits_{{\mathbb R}^n} \,|D u|^2\,dx\ge$$$$\ge -(\gamma-1)\,
E_i(t) \, f_1(t)+ K \,(E_i(t))^{-\frac{n-2}{n(\gamma-1)}}\ge -L_1
\,t^\beta \,f_1(t) + L_2 \, t^{\frac{n-2}{n}}, $$

\medskip

$$
E_i'(t)\ge -\frac{\gamma-1}{\mu+\lambda} E_i(t) f_2(t)+ \mu
\int\limits_{{\mathbb R}^n} \,|\nabla u|^2\,dx\ge$$$$
 \ge -\frac{\gamma-1}{\mu+\lambda} E_i(t) f_2(t)+ K (E_i(t))^{-\frac{n-2}{n(\gamma-1)}}\ge
  -L_3 \,t^\beta\,f_2(t) + L_2 \, t^{\frac{n-2}{n}},$$
with positive constants $L_1, L_2, L_3$ where $\beta=-n(\gamma-1),$
for $\gamma\le 1+\frac{2}{n},$ and $\beta = -2$, otherwise. If the
growth rate of $f_1(t)$ and $f_2(t)$ is less then prescribed in
Theorem 3.1 statement, we get a contradiction. Thus, the theorem is
proved.

\begin{remark} As follows from (3.22) and Lemmas 3.2, 3.3, 3.4,
the requirement of incompressibility ${\rm div}_x u=0$ signifies
that in conditions of Theorem 2.2 the solution to (NS*) loses its
initial smoothness within a finite time.

\end{remark}

\subsection{Long-time behavior of solution and the blow up}

In \cite{ChoJa} it was found that there exists no global smooth
solution to (NS) such that
$$\limsup\limits_{t\to\infty}\,\left\|t\,\frac{(u,x)}{1+|x|^2} \right \|_{L^\infty({\mathbb
R}^n)}<1.\eqno(3.23)
$$
One can also derive that there exists no global smooth solution to
(NS) such that
$$
\limsup\limits_{t\to\infty}\,\left\| u \right \|_{L^\infty({\mathbb
R}^n)}<\frac{|P|}{m}= const\eqno(3.24)
$$
or
$$
\limsup\limits_{t\to\infty}\,\left(\int\limits_{t_0}^t\,\left\|{\rm
div} u \right \|_{L^\infty({\mathbb R}^n)}(\tau)\, d\tau - {n}\,{\ln
t}\right) < const, \quad t_0\ge 0.\eqno(3.25)
$$
Condition (3.24) follows from the H\"older inequality:
$$
|P|=\Big|\int\limits_{{\mathbb R}^n}\,\rho u\,dx\Big|\le m\,\left\|
u \right \|_{L^\infty({\mathbb R}^n)}.
$$

To prove (3.25) (e.g. for $\gamma\le 1+\frac{2}{n}$) we note that
(3.22) results in
$$
E_i'(t)\ge -\int\limits_{{\mathbb R}^n}\,p\, \, {\rm div} u\,dx\ge
-(\gamma-1)\, E_i(t)\,\left\|{\rm div} u \right
\|_{L^\infty({\mathbb R}^n)}. $$ Further, integrating from any
$t_0\ge 0$ gives
$$
\ln E_i(t)\ge -(\gamma-1)\, \int\limits_{t_0}^t\,\left\|{\rm div} u
\right \|_{L^\infty({\mathbb R}^n)}(\tau)\, d\tau+\ln
E_i(t_0).\eqno(3.26)
$$
Further, inequalities (3.3), (3.9) give the following estimate for
sufficiently large $t$:
$$
E_i(t)\le C_2\,(\mathcal E
t^2+F(0)t+G(0))^{-(\gamma-1)n/2}.\eqno(3.27)
$$
Estimate (3.25) follows from (3.26) and (3.27) immediately.

However,  condition (3.23) just as  conditions (3.24) and (3.25)
does not use the fact that the velocity  belongs to the space
$H^1({\mathbb R}^n).$ These conditions can be applied for the case
of zero coefficients of viscosity, i.e. for the gas dynamics
equations. Moreover, it is possible to construct global in time
exact solutions to the gas dynamics equations with the velocity that
increases by modulus as $|x|\to\infty$ and has the  form $u=A(t) \,
x,\quad A(t)\sim t^{-1},\, t\to\infty\,$ (see e.g.
\cite{RozNova},\cite{Ro} for details) such that conditions (3.23),
(3.24), (3.25) become equalities. These solutions satisfy (NS) as
well.  As follows from (3.25), for  smooth solutions to (NS) the
function $\int\limits_{t_0}^t\,\left\|{\rm div} u \right
\|_{L^\infty({\mathbb R}^n)}(\tau)\, d\tau \ge {n}\,{\ln t}+
const,\, t_0\ge 0,\,$ and the comparison with with the statement of
Theorem 2.2 shows that for solutions of class $H^1({\mathbb R}^n)$
this estimate is very far to be exact.

The author thanks Profs. S.Albeverio and A.A.Zlotnik for a helpful
discussion.

\end{document}